\documentclass[12pt]{article}
\usepackage{times}

\usepackage{amssymb}
\newcommand{\nc}[2]{\newcommand{#1}{#2}}
\nc{\lra}{\longrightarrow}
\nc{\ra}{\rightarrow}
\nc{\ci}{\circ}
\nc{\beq}{\begin{equation}}
\nc{\eeq}{\end{equation}}
\nc{\bea}{\begin{eqnarray*}}
\nc{\eea}{\end{eqnarray*}}
\newtheorem{definition}{Definition}
\newtheorem{proposition}{Proposition}
\newtheorem{theorem}{Theorem}
\newtheorem{lemma}{Lemma}
\nc{\ble}{\begin{lemma}}
\nc{\ele}{\end{lemma}}
\nc{\bde}{\begin{definition}}
\nc{\ede}{\end{definition}}
\nc{\bpr}{\begin{proposition}}
\nc{\epr}{\end{proposition}}
\nc{\bth}{\begin{theorem}}
\nc{\ethe}{\end{theorem}}
\nc{\ot}{\otimes}
\nc{\Hom}{{\rm Hom}}
\def\<{\langle}
\def\>{\rangle}
\nc{\ba}{\begin{array}}
\nc{\ea}{\end{array}}
\def\Z{{\mathbb Z}}

\def\C{{\mathbb C}}
\def\T{{\mathbb T}}

\def\cH{{\cal H}}

\def\cO{{\cal O}}
\def\cT{{\cal T}}
\def\id{{\rm id}}
\begin{document}
\title{A locally trivial quantum Hopf bundle}
\author{R. Matthes
\\
Fachbereich Physik der TU Clausthal\\ Leibnizstr. 10,
D-38678 Clausthal-Zellerfeld, Germany}
\date{}
\maketitle 
\begin{abstract}
We describe a locally trivial quantum principal $U(1)$-bundle over the
quantum space $S^2_{pq}$ which is a noncommutative analogue of the
usual Hopf bundle. We also provide results concerning the structure of
its total space algebra (irreducible $*$-representations and topological
$K$-groups) and its Galois aspects (Galois property, existence of a strong 
connection, non-cleftness).
\end{abstract}

\section{Introduction}

In this note, we describe an example of a principal bundle
in the setting of noncommutative geometry, which meets 
two possible (still provisional) definitions:
It is a locally trivial quantum principal bundle in the sense of 
\cite{bk96} as well as a Hopf-Galois extension \cite{m-s93}.
Besides giving the definition of these notions and a description 
of the bundle \cite{cm00}, \cite{cm02},
we provide a list of results obtained in \cite{hms} 
concerning the structure of the total space algebra
and the Galois aspects of the bundle. 
\section{Quantum principal bundles}
\subsection{Hopf-Galois extensions}
\label{hg}
Dualizing the corresponding classical structure ``\`a la Gelfand-Neumark'',
one arrives at the following items which show up
in the definition of quantum principal bundles:
\begin{itemize}
\item
There is some algebra $P$ replacing the total space of a principal bundle.
\item
There is some Hopf algebra $H$ replacing the structure group, coacting on $P$ 
on the right,
i.e., there is an algebra homomorphism $\Delta_R:P\ra P\ot H$
with $(\Delta_R\ot\id)\ci\Delta_R=(\id\ot\Delta)\ci\Delta_R$
and $(\id\ot\varepsilon)\ci\Delta_R=\id$.
\item
There is another algebra replacing the base space, which coincides with
the subalgebra of coinvariants of the coaction of $H$ on $P$, 
$B=P^{coH}:=\{p\in P\:|\:\Delta_R(p)=p\ot 1\}$.
The bundle projection is the 
embedding $B\subset P$, denoted by $\iota:B\ra P$.
\end{itemize}
$B\subset P$ is called $H$-extension in the above context \cite{m-s93}.
For a classical principal bundle with base space $M$, total space $P$ and
structure group $G$, the right action is assumed to be free.
This assumption can be restated as bijectivity of the map
$X\times G\ra X\times_MX,~~(x,g)\mapsto (x,xg)$. At the level of algebras,
this means bijectivity of the map
\[
can: P\ot_BP\lra P\ot H,~~p\ot p'\mapsto pp'_{(0)}\ot p'_{(1)}.
\]
Here we use Sweedler notation, $\Delta_R(p)=p_{(0)}\ot p_{(1)}$. 
An $H$-extension
is called Hopf-Galois if $can$
is bijective.
This is
essentially the notion of an algebraic quantum principal bundle (see,
e.g., \cite{bm93}). 
\subsection{Locally trivial quantum principal bundles}\label{lt}
There is another approach to quantum principal
bundles emphasizing the idea of gluing which is behind the definition of 
classical fibre bundles \cite{bk96}, \cite{cm02}. 
In order to state this definition, we need an
algebraic notion of covering:
 
A covering of an algebra $B$ is a family $(J_i)_{i\in I}$ of ideals with zero
intersection. Let $\pi_i:B\ra B_i:=B/J_i$, $\pi^i_j:B_i\ra B_{ij}:=
B/(J_i+J_j)$ be the quotient maps.
A covering $(J_i)_{i\in I}$ is called complete if the homomorphism
$
B\ni b\mapsto (\pi_i(b))_{i\in I}\in
\{(b_i)_{i\in I}\in\prod_{i\in I}B_i\:|\:\pi^i_j(b_i)=\pi^j_i(b_j)\}
$
is surjective (it is always injective). 
Finite coverings by closed ideals in C*-algebras and two-element coverings 
are always complete. 
A locally trivial $H$-extension is an
$H$-extension $B\subset P$ supplied
with the following local data:\\
(i) $B$ has a complete finite covering $(J_i)_{i\in I}$.\\
(ii) There are given surjective homomorphisms
$\chi_i:P\ra B_i\ot H$ (local trivializations) such that\\
\hspace*{1cm}(a) $\chi_i\ci\iota=\pi_i\ot 1$ ($\iota:B\ra P$),\\
\hspace*{1cm}(b) $(\chi_i\ot\id)\ci\Delta_R=(\id\ot\Delta)\ci\chi_i$
(right colinearity),\\
\hspace*{1cm}(c) $(\ker\chi_i)_{i\in I}$ is a complete covering of $P$.
\vspace{.2cm}

\noindent
As in the classical situation, locally trivial bundles can be
reconstructed from transition functions related
to the covering of the base algebra. More precisely, every locally trivial
principal fibre bundle with fixed base algebra $B$ and Hopf algebra $H$ is 
determined by the following data:
\begin{itemize}
\item
a complete finite covering
$(J_i)_{i\in I}$
\item
a family  of transition functions, i.e., of homomorphisms 
$\tau_{ij}:H\ra Z(B_{ij})$ (center) fulfilling
$\tau_{ii}=1\varepsilon$, $
\tau_{ji}\ci S=\tau_{ij}$ ($S$ the antipode of $H$),
and the cocycle condition 
$\pi^{ij}_k\ci\tau_{ij}=m_{B_{ijk}}\ci((\pi^{ik}_j\ci\tau_{ik})
\ot(\pi^{jk}_i\ci\tau_{kj}))\ci\Delta$.
\end{itemize}
The total space algebra is then given as the gluing
\[
P=\{(f_i)_{i\in I}\in\oplus_{i\in I}B_i\ot H\:|\:(\pi^i_j\ot\id)(f_i)=
\varphi_{ij}\ci(\pi^j_i\ot\id)(f_j)\},
\]
where $\varphi_{ij}(b\ot h)=b\tau_{ji}(h_{(1)})\ot h_{(2)}$.
The remaining data of the corresponding locally trivial $H$-extension
are as follows:
\[
\Delta_R((f_i)_{i\in I})=((\id\ot\Delta)(f_i))_{i\in I},~~
\chi_i((f_i)_{i\in I})=f_i,~~
\iota(b)=(\pi_i(b)\ot 1)_{i\in I}.
\]
\section{Description of
the locally trivial $\boldmath{U(1)}$-bundle 
$\boldmath{S^3_{pq}\rightarrow S^2_{pq}}$}
\subsection{Quantum discs}
We use the following subfamily of a two-parameter family of quantum
discs defined in \cite{kl93} whose $*$-algebra is
$
\cO(D_q):= \C\langle x,x^*\rangle/(x^*x-qxx^*-(1-q)),~~0<q<1.
$
The irreducible
$*$-representations of $\cO(D_q)$ are an $S^1$-family of one-dimensional 
representations,
given by
$\pi_\theta(x)=e^{i\theta}$ (classical points),
and an infinite-dimensional representation $\pi_q$ in a separable Hilbert space
representing the generator $x$ as a one-sided weighted shift.
The classical points define an embedding of $S^1$ into $D_q$,
i.e.,
$
\cO(D_q)\ni x\stackrel{\phi_q}{\lra}\underline{u}\in
\cO(S^1):=\C\langle \underline{u},\underline{u}^*\rangle/(\underline{u}^*
\underline{u}-1,\underline{u}\hspace{1mm}\underline{u}^*-1).
$
\noindent
Since $\|\pi(x)\|=1$ for any $*$-representation of $\pi$ in some $B(\cH)$, the
$C^*$-closure $C(D_q)$ of $\cO(D_q)$ is well-defined (using bounded 
$*$-representations).
One knows that $C(D_q)\simeq\cT$ (Toeplitz or shift algebra).
Using the above-mentioned irreducible $*$-representations, one may heuristically
interprete $D_q$ as
a diffuse membrane spanned by a classical $S^1$.
\subsection{Quantum two-spheres (quantum cones)}
They are defined as a gluing of two quantum discs along the classical
``boundary'' $S^1$:
$
\cO(S^2_{pq}):= \cO(D_p)\oplus_\phi\cO(D_q)=\{(f,g)\in
\cO(D_p)\oplus\cO(D_q)\:|\:\phi_p(f)=\phi_q(g)\},
~~0<p,q<1.
$
The $*$-algebra
$
\cO(S^2_{pq})$ can be identified with the quotient of the free algebra
generated by $f_1,f_1^*,f_0$
by the ideal $J$ defined by the relations
$
f_0^*=f_0,~~
f_1^*f_1-qf_1f_1^*=(p-q)f_0+(1-p)1,~~
(1-f_0)(f_1f_1^*-f_0)=0.
$
There are an $S^1$-family of one-dimensional 
and two nonequivalent infinite dimensional $*$-representations
in a separable Hilbert space. The latter represent $f_0$ as a 
diagonal operator
and $f_1$ as a one-sided weighted shift. 
Again $\|\rho(f_0)\|=\|\rho(f_1)\|=1$ for any bounded $*$-representation.
The $C^*$-closure $C(S^2_{pq})$
is defined using such representations.
One knows
$
C(S^2_{pq})\simeq C(D_p)\oplus_\phi C(D_q),
\simeq C(S^2_{\mu c}),~~|\mu|<1,c>0
$
(Podle\'s spheres \cite{p-p87}). Thus, the glued two-spheres are homeomorphic to the
so-called equilateral Podle\'s spheres.
Using, as for the disc, the irreducible representations, one may visualize
$S^2_{pq}$ as a top of a diffuse 
cone, with edge $S^1$.

\subsection{The $\cO(U(1))$-extension $\cO(S^2_{pq})\subset \cO(S^3_{pq})$}
Note that $\cO(S^2_{pq})= \cO(D_p)\oplus_\phi\cO(D_q)$ has a canonical covering
consisting of the kernels of the first and second projections,
$J_1=\ker pr_1$, $J_2=\ker pr_2$. One has canonical identifications
$\cO(S^2_{pq})/J_1= \cO(D_p)$, $\cO(S^2_{pq})/J_2= \cO(D_q)$,
$\cO(S^2_{pq})/(J_1+J_2)= \cO(S^1)$. These are the data of the base algebra.
The desired extension results from gluing $\cO(D_p)\ot \cO(U(1))$ and $\cO(D_q)\ot
\cO(U(1))$ by means of one transition function
$
\tau:\cO(U(1))\lra\cO(S^1),~~ u\mapsto\underline{u},
$
following the general method of Subsection \ref{lt}.
The corresponding gluing $\cO(S^3_{pq})$ of two quantum solid tori along
their set $\T^2$
of classical points is fully analogous to the geometrical picture in the case
of the usual $U(1)$-Hopf bundle (Heegard splitting of $S^3$). It turns out that $\cO(S^3_{pq})$
is isomorphic to the quotient of the free $*$-algebra generated by
$a, b$ by the ideal generated by the relations
\[
ab=ba,~~ab^*=b^*a,~~a^*b^*=b^*a^*,~~a^*b=ba^*,
\]
\[
a^*a-qaa^*=1-q,~~~~b^*b-pbb^*=1-p,
\]
\[
(1-aa^*)(1-bb^*)=0.
\]
The structural $*$-homomorphisms of the locally trivial $U(1)$-extension
 in terms of the
generators $a,b$ are:
\[
\Delta_R(a)=a\ot u,~~~\Delta_R(b)=b\ot u^*,
\]
\[
\chi_p(a)=1\ot u, ~~~\chi_p(b)=x\ot u^*,~~~
\chi_q(a)=y\ot u,~~~\chi_q(b)=1\ot u^*,\]
\[
\iota(f_1)=ba,~~~\iota(f_0)=bb^*.
\]

\section{Further results}
\subsection{Structure of $S^3_{pq}$}
\begin{itemize}
\item
The classes of irreducible $*$-representations of $\cO(S^3_{pq})$
in bounded operators are 
classified:
There is a $\T^2$-family of one dimensional representations and two 
$S^1$-families of infinite-dimensional representations in a separable
Hilbert space. In the first of these two families, $a$ is a multiple of the
unit operator, and $b$ is a one-sided weighted shift. In the second
family $a$ and $b$ exchange their roles.
Since again the norms of $a$ and $b$ are 1 in any bounded representation, one
can define the $C^*$-algebra $C(S^3_{pq})$ using such representations.
\item
A vector space basis of $\cO(S^3_{pq})$ can be exhibited.
\item
$C(S^3_{pq})$ is a 2-graph $C^*$-algebra.
\item
The $K$-groups of $C(S^3_{pq})$ coincide with the
$K$-groups of the classical $S^3$, i.e.,\\
$K_0(C(S^3_{pq}))=K_1(C(S^3_{pq}))=\Z$.
\end{itemize}
\subsection{Hopf-Galois (bundle) aspects}
\begin{itemize}
\item
The $\cO(U(1))$-extension $\cO(S^2_{pq})\subset\cO(S^3_{pq})$ has the Galois
property. (Idea of proof: Find a lift $l$ of the translation map and use a
general argument of Schneider.)
\item
The lift $l$ of the translation map is a strong connection in the sense of
\cite{h-pm96}. Consequently, the $\cO(U(1))$-extension 
$\cO(S^2_{pq})\subset\cO(S^3_{pq})$ is relatively projective \cite{bh}.
\item
As a further consequence
of the existence of a strong connection, all associated modules 
(vector bundles) are finitely generated projective. In particular,
using the strong connection one can for any winding number give explicitely 
a projector matrix
corresponding to the associated line bundle.
\item
The $\cO(U(1))$-extension
$\cO(S^2_{pq})\subset\cO(S^3_{pq})$ is non-cleft (not a crossed product).
This is proved using a trace on $\cO(S^2_{pq})$, which is
defined as the operator trace composed with the difference of the two
irreducible infinite dimensional representations. The Chern-Connes
pairing of this trace with the $K_0$-class of the projector
defining the associated line bundle with winding number -1 just gives this 
number, which proves
the above claim (cf. \cite{hm99}). 
\end{itemize}
{\bf Acknowledgements}:
This work was supported by the Deutsche Forschungsgemeinschaft 
and the Mathematisches Forschungsinstitut Oberwolfach, where this note was
completed during a stay under
the Research in Pairs programme. Also, it is a pleasure to thank D. Calow,
P.M. Hajac and W. Szymanski for many hours of discussion and joint work.

\end{document}